\documentclass[graybox]{svmult}

\usepackage[pagewise]{lineno}   

\usepackage{mathptmx}       
\usepackage{helvet}         
\usepackage{courier}        
\usepackage{type1cm}        
\usepackage{makeidx}         
\usepackage{graphicx}        
\usepackage{multicol}        
\usepackage[bottom]{footmisc}

\usepackage{pstricks-add}
\usepackage{verbatim}
\usepackage{amssymb}
\usepackage{latexsym}
\usepackage[all]{xy}

\renewcommand{\geq}{\geqslant}

\newcommand{\K}{\mathcal{K}}

\newrgbcolor{fffftt}{1 1 0.2}
\newrgbcolor{ffttqq}{1 0.2 0}
\newrgbcolor{ffffzz}{1 1 0.6}
\newrgbcolor{qqcctt}{0 0.8 0.2}

\makeindex             

\begin{document}


\title*{Polygonal Complexes and Graphs for Crystallographic Groups}

\author{Daniel Pellicer and Egon Schulte}
\institute{Daniel Pellicer, \at Centro de Ciencias Matematicas, National University of Mexico, CP 58087, Morelia, Michoacan, Mexico,
\email{pellicer@matmor.unam.mx}
\and
Egon Schulte, \at Northeastern University, Department of Mathematics, Boston, MA 02115, USA, \email{schulte@neu.edu}}
%
%
\maketitle

\abstract{The paper surveys highlights of the ongoing program to classify discrete polyhedral structures in Euclidean $3$-space by distinguished transitivity properties of their symmetry groups, focussing in particular on various aspects of the classification of regular polygonal complexes, chiral polyhedra, and more generally, two-orbit polyhedra.}


\section{Introduction}
\label{intro}

The study of highly symmetric discrete structures in ordinary Euclidean 3-space $\mathbb{E}^3$ has a long and fascinating history tracing back to the early days of geometry. With the passage of time, various notions of discrete structures with properties similar to those convex polyhedra have attracted attention and have brought to light new exciting figures intimately related to finite or infinite groups of isometries.

A radically new ``skeletal" approach to polyhedra in $\mathbb{E}^3$ was pioneered by Gr\"unbaum~\cite{gr1} in the 1970's, building on Coxeter's work~\cite{coxeter}. A polyhedron is viewed as a finite or infinite periodic geometric (edge) graph in space equipped with additional structure imposed by the faces, and its symmetry is measured by transitivity properties of its geometric symmetry group. For example, the geometric graph of the cube carries four {\em Petrie polygons\/}, that is, polygons for which any two, but no three, consecutive edges belong to the same square of the cube. The geometric graph of the cube with its four hexagonal Petrie polygons constitutes one of the new regular polyhedra introduced by Gr\"unbaum. Throughout this paper we shall adopt this notion of polyhedron. 

Since the mid 1970's, there has been a lot of activity in this area, beginning with the full enumeration of the ``new" regular polyhedra by Gr\"unbaum~\cite{gr1} and Dress~\cite{d1,d2} by around 1980 (see also McMullen-Schulte~
\cite[Ch.~7E]{arp} or \cite{ordinary} for a faster method for arriving at the complete list); moving on to the full enumeration of the chiral polyhedra in \cite{chiral1,chiral2} by around 2005; and continuing with the enumeration of certain classes of regular polyhedra and polytopes in higher-dimensional spaces by McMullen~\cite{pm,grp}.

While all these structures have the essential characteristics of polyhedra and polytopes, the more general class of discrete ``polygonal complexes" in $3$-space is a hybrid of polyhedra and incidence geometries (see \cite{bu2}). Every edge of a polyhedron belongs to precisely two edges, whereas an edge of a polygonal complex is surrounded by any number at least two. For example, the geometric edge graph of the cube endowed with the six squares and four Petrie polygons as faces constitutes a polygonal complex where every edge belongs to precisely four faces. In very recent joint work, we obtained a complete enumeration of the regular polygonal complexes in $\mathbb{E}^3$ (see \cite{pelsch1,pelsch2}). These are periodic structures with crystallographic symmetry groups exhibiting interesting geometric, combinatorial, and algebraic properties.

The purpose of this paper is to exhibit some of the highlights of the ongoing program to classify discrete structures built from vertices, edges and faces in Euclidean $3$-space according to transitivity properties of their symmetry groups. We center our attention on the recent classification of regular polygonal complexes, chiral polyhedra, and more generally, two-orbit polyhedra.

In Sections~\ref{termin} and \ref{symgroup}, we review basic terminology about polygonal complexes and describe structure results for the symmetry group of regular polygonal complexes. This is followed, in Section~\ref{regpol}, by a brief description of the complete enumeration of regular polyhedra, seen from the perspective of regular polygonal complexes. Then Sections~\ref{nonsim} and \ref{simco} give an account of the regular polygonal complexes which are not polyhedra. In the last two sections we study certain kinds of two-orbit polyhedra in~$\mathbb{E}^3$, beginning with a review of the enumeration of chiral polyhedra. Finally, Section~\ref{2orbpoly} briefly summarizes the recent classification of regular polyhedra of index $2$, obtained in Cutler~\cite{cut} and \cite{cutsch}; these form a distinguished class of two-orbit polyhedra in $\mathbb{E}^3$.

\section{Some terminology}
\label{termin}

Informally, a polygonal complex is a discrete structure in $\mathbb{E}^3$ consisting of vertices (points), joined by edges (line segments) assembled in careful fashion into faces (polygons, allowed to be finite or infinite), with at least two faces on each edge. Our concepts of face (polygon) and polyhedron generalize those of convex polygon and convex polyhedron.

For our purposes, a {\em finite polygon\/} $(v_1, v_2, \dots, v_n)$ in $\mathbb{E}^3$ is a figure consisting of $n$ distinct points $v_1, \dots, v_n$, together with the line segments $(v_i, v_{i+1})$, for $i = 1, \dots, n-1$, and $(v_n, v_1)$. Similarly, an {\em infinite polygon\/} is a figure made up from an infinite sequence of distinct points $(\dots, v_{-2},v_{-1}, v_0, v_1, v_2,\dots)$, and line segments $(v_i, v_{i+1})$ for each $i$, such that each compact subset in $\mathbb{E}^3$ meets only finitely many line segments. In either case the points and line segments are referred to as the {\em vertices\/} and {\em edges\/} of the polygon, respectively.

By a {\em regular} polygon we mean a finite or infinite polygon such that its geometric symmetry group, restricted to the affine hull of the vertices, is a finite or infinite dihedral group acting transitively on the set of incident vertex-edge pairs, called {\em flags}. This definition covers not only the traditional (planar) convex regular polygons but also permits star-polygons, skew polygons, zigzags, or helices as regular polygons. A {\em star polygon} has the same vertices as a convex regular polygon; its edges connect vertices of the convex regular polygon that are a fixed number of steps apart on the boundary. A {\em skew polygon} lives properly in $\mathbb{E}^3$ and can be obtained from a planar finite (convex or star-) polygon by raising every other vertex  perpendicularly by the same amount (thus doubly covering the original polygon if the number of vertices was odd); the vertex set then is contained in two parallel planes and every edge goes from a vertex in one plane to a vertex in the other plane. A {\em linear apeirogon} is an infinite polygon obtained by tessellating a line with line segments (usually of the same size). Linear apeirogons will not occur as faces of the geometric objects described in this paper, since no non-trivial connected structure can be assembled only from linear building blocks. A {\em zigzag} is a planar infinite polygon obtained from a linear apeirogon in a similar way as a skew polygon is obtained from a planar finite polygon; its vertices lie on two parallel lines, and its edges connect vertices on different lines. Finally, a {\em helix} is an infinite non-planar polygon and it can be thought as a spring rising above a finite planar (convex or star) polygon; more precisely, the orthogonal projection onto its axis gives a linear apeirogon, and the orthogonal projection along its axis gives a finite planar (convex or star-) polygon.

A {\em polygonal complex}, or simply {\em complex}, $\K$ in $\mathbb{E}^3$ consists of a set $\mathcal{V}$ of points, called {\em vertices}, a set $\mathcal{E}$ of line segments, called {\em edges}, and a set $\mathcal{F}$ of polygons, called {\em faces}, such that the following properties are satisfied. The graph defined by $\mathcal{V}$ and $\mathcal{E}$, called the {\em edge graph\/} of $\K$, is connected. Moreover, the vertex-figure of $\K$ at each vertex of $\K$ is connected. Here the {\em vertex-figure\/} of $\K$ at a vertex $v$ is the graph, possibly with multiple edges, whose vertices are the neighbors of $v$ in the edge graph of $\K$ and whose edges are the line segments $(u,w)$, where $(u, v)$ and $(v, w)$ are edges of a common face of $\K$. (Note that this is a small change over \cite{pelsch1,pelsch2}, where a complex was required to have exactly $r$ faces on each edge, for some fixed number $r \geq 2$. However, for regular complexes with at least two faces meeting at an edge the two definitions are equivalent.) All polygonal complexes studied in this paper have at least two faces on each edge. Finally, $\K$ is {\em discrete\/}, in the sense that each compact subset of $\mathbb{E}$ meets only finitely many faces of $\K$. A complex with exactly two faces on each edge (that is, $r=2$) is also called a {\em polyhedron\/}. Note that this definition extends the notion of convex polyhedron, where the faces are convex (finite and planar) and the polyhedron itself is also finite. Polyhedra (finite or infinite) with high symmetry properties have been extensively studied in \cite[Ch.~7E]{arp} and \cite{gr1,ordinary,pelwei,chiral1,chiral2}.

A polygonal complex $\K$ is said to be {\em regular\/} if its geometric symmetry group $G:=G(\K)$ is transitive on the incident vertex-edge-face triples, called {\em flags}. The faces of a regular complex are necessarily regular polygons.
The vertex-figures are finite (flag-transitive) graphs with single or double edges; double edges occur precisely when any two adjacent edges of a face of $\K$ are adjacent edges of another (then uniquely determined) face of $\K$. If $\K$ is not a polyhedron, then $G$ is infinite and affinely irreducible, that is, $G$ is a standard crystallographic group (see \cite{pelsch1}). In particular, there are no finite regular complexes other than polyhedra. The Platonic solids are the most natural examples of regular polyhedra, and the $2$-skeleton of the tessellation of $\mathbb{E}^3$ by cubes is the most natural example of a regular polygonal complex which is not a polyhedron.

Regular polygonal complexes in $\mathbb{E}^3$ can be viewed as $3$-dimensional (discrete faithful) Euclidean realizations of regular incidence complexes of rank $3$ with polygonal faces (see \cite{kom1,kom2}). Our description of the symmetry groups will exploit this fact. In particular, the regular polyhedra in $\mathbb{E}^3$ are precisely the $3$-dimensional discrete faithful Euclidean realizations of abstract regular polyhedra (abstract regular $3$-polytopes); for more details, see \cite[Ch. 7E]{arp} and \cite{ms3}.

Every regular polyhedron has the property that all its faces have the same number $p$ of edges, and all its vertices have the same degree $q$. Polyhedra with this property are called {\em equivelar}, and their {\em Schl\"afli type} (or {\em Schl\"afli symbol}) is defined to be $\{p,q\}$. When the faces of an equivelar polyhedron are zigzags or helices, the first entry $p$ is $\infty$; however, since we only consider discrete structures, $q$ is always finite. Similarly, in the {\em Schl\"afli symbol\/} $\{p,q,r\}$ of a regular rank 4 polytope (a combinatorial structure contructed from vertices, edges, polygons and polyhedra) the first two entries give the Schl\"afli type $\{p,q\}$ of any of its rank 3 faces, while the last entry $r$ is the number of rank 3 faces meeting around each edge (so that the vertex-figures have Schl\"afli type $\{q,r\}$).

In later sections we also meet various kinds of less symmetric polygonal complexes (in fact, polyhedra) in $\mathbb{E}^3$. These have more than one flag orbit under the symmetry group. A particularly interesting case arises when there are just two flag orbits. We say that a polygonal complex $\K$ is a {\em $2$-orbit polygonal complex\/} if $K$ has precisely two flag-orbits under $G$; in this case, if $\K$ is also a polyhedron, we call $\K$ a {\em $2$-orbit polyhedron\/}. The cuboctahedron and icosidodecahedron are simple examples of two-orbit polyhedra.

There are different kinds of $2$-orbit polyhedra in $\mathbb{E}^3$. Recall that two flags of a polyhedron $\K$ are called {\em $i$-adjacent\/}, with $i=0$, $1$, or $2$ respectively, if they differ precisely in their vertices, edges, or faces (see \cite[Ch. 2]{arp}). Thus, two flags are $1$-adjacent if they have the same vertices and same faces, but different edges. Note that flags of polyhedra have unique $i$-adjacent flags for each $i$; for polygonal complexes which are not polyhedra, this still is true for $i=0,1$ but not for $i=2$. Now $2$-orbit polyhedra naturally fall into different classes indexed by proper subsets $I$ of $\{0,1,2\}$ (see Hubard~\cite{hub} and \cite{hubsch}). In particular, a $2$-orbit polyhedron $\K$ is said to belong to the {\em class\/} $2_I$ if $I$ consists precisely of those indices $i$ such that any two $i$-adjacent flags lie in the same flag-orbit under $G$. The cuboctahedron and the icosidodecahedron are examples of two orbit polyhedra in class $2_{\{0,1\}}$. When $I=\emptyset$ this gives the class $2_{\emptyset}$ of {\em chiral\/} polyhedra. Thus a polyhedron $\K$ is chiral if and only if $\K$ has two flag orbits under $G$ such that any two adjacent flags lie in distinct orbits. (The case $I=\{0,1,2\}$ is excluded here, as it describes the regular polyhedra.)

\section{The symmetry group}
\label{symgroup}

The symmetry group $G=G(\K)$ of a regular complex $\K$ in $\mathbb{E}^3$ either acts regularly on the set of flags or has flag-stabilizers of order $2$. We call $\K$ {\em simply flag-transitive\/} if its (full) symmetry group $G$ acts regularly on the flags of $\K$; in other words, $G$ is simply transitive on the flags of $\K$. Note that a regular complex that is not simply flag-transitive can (in fact, always does) have a subgroup (of index $2$) that acts simply flag-transitively. Each regular polyhedron, finite or infinite, is a simply flag-transitive regular polygonal complex.

The group $G$ always has a well-behaved system of generators or generating subgroups, regardless of whether $\K$ is simply flag-transitive or not. Suppose $\Phi := \{F_0, F_1, F_2\}$ is a fixed, or {\em base}, flag of $\K$, consisting of a vertex $F_0$, an edge $F_1$, and a face $F_2$. For each $\Psi\subseteq\Phi$ we let $G_{\Psi}$ denote the stabilizer of $\Psi$ in $G$.  Moreover, for $i=0,1,2$ we set $G_i := G_{\{F_j, F_k\}}$, where $i,j,k$ are distinct, and write $G_{F_i}:=G_{\{F_i\}}$ for the stabilizer of $F_i$ in $G$. Then $G_\Phi$ is the stabilizer of $\Phi$ and has order $1$ or $2$; in particular,
\[G_\Phi = G_0 \cap G_1 = G_0 \cap G_2 = G_1 \cap G_2.\]
The stabilizers $G_0,G_1,G_2$ form a generating set of subgroups for $G$, with the property that $G_0 \cdot G_2 = G_2 \cdot G_0 = G_{F_1}$ and $G_\Psi = \langle G_j \,|\, F_j \notin \Psi\rangle$ for each $\Psi \subseteq \Phi$. Moreover,
\[ \langle G_j \,|\, j \in I\rangle \cap \langle G_j \,|\, j\in J\rangle =
\langle G_j \,|\, j \in I\cap J\rangle \quad (I, J \subseteq \{0, 1, 2\}).\]
These statements about generating subgroups of $G$ are particular instances of similar such statements about flag-transitive subgroups of automorphism groups of regular incidence complexes of rank $3$ (or higher) obtained in \cite[\S2]{kom2} (and also described in \cite[pp. 33,34]{arp} for polyhedra).

From the base vertex $F_0$ and the symmetry group $G$ of a regular complex $\K$, with generating subgroups $G_0$, $G_1$, $G_2$, we can reconstruct $\K$ by the following procedure, often called {\em Wythoff's construction}. First  observe that the base edge $F_1$ of $\K$ is determined by the pair of vertices $\{F_0, F_0 G_0\}$. Similarly, the vertex- and edge-sets, respectively, of the base face $F_2$ of $\K$ are just $\{F_0 S \mid S\in \langle G_0, G_1\rangle\}$ and $\{F_1 S \mid S \in \langle G_0,G_1\rangle\}$. This recovers the base flag of $\K$. Finally, the set of $i$-faces of $\K$ is just $\{F_i S \mid S \in G\}$ for each $i=0,1,2$.

Most regular complexes $\K$ in $\mathbb{E}^3$ are infinite and have an affinely irreducible infinite discrete group of isometries as a symmetry group. In this case $G$ is a crystallographic group (that is, $G$ admits a compact fundamental domain). Then the
Bieberbach theorems tell us that $G$ contains a translation subgroup (of rank $3$) such that the quotient of $G$ by this subgroup is finite (see \cite[\S7.4]{ratcliffe}). If $R: x \mapsto xR' + t$ is a general element of $G$, with $R'$ in $\rm{O}(3)$, the orthogonal group of $\mathbb{E}^3$, and $t$ a translation vector in $\mathbb{E}^3$ (that we may also view as a translation), then the mappings $R'$ clearly form a subgroup $G_*$ of $\rm{O}(3)$, called the {\em special group\/} of $G$. Now if $T(G)$ denotes the full translation subgroup of $G$ (consisting of all translations in $G$), then
\[ G_*=G/T(G),\]
so in particular, $G_*$ is a finite group. Thus $G_*$ is among the finite subgroups of $\rm{O}(3)$, which are known (see \cite{grove}). The special group of any irreducible infinite discrete group of isometries in $\mathbb{E}^2$ or $\mathbb{E}^3$ never contains rotations of periods other than $2, 3, 4$, or $6$, and period $6$ only occurs for $\mathbb{E}^2$ (see \cite[p.220]{arp} and \cite[Lemma 3.1]{chiral2}).

The full translation subgroup of the symmetry group $G$ of a regular complex $\K$ (and often the vertex set of $\K$ itself) is given by a $3$-dimensional lattice in $\mathbb{E}^3$. We frequently meet the lattices $\Lambda_{\bf a}$ that are generated by a single vector ${\bf a} : = (a^k,0^{3-k})$ and its images under permutations and changes of sign of coordinates; here $a>0$ and $k=1,2,3$ (and ${\bf a}$ has $k$ entries $a$ and $3-k$ entries $0$). When $a=1$ and $k=1$, $2$ or $3$, respectively, these are
the standard {\em cubic lattice\/} $\mathbb{Z}^{3}$, the {\em face-centered cubic lattice\/}, and the {\em body-centered cubic\/} lattice.

\section{Regular polyhedra}
\label{regpol}

The regular polyhedra in space are also known as the {\it Gr\"unbaum-Dress polyhedra\/} (see \cite{sympopo}). It is convenient to separate them from the simply flag-transitive regular complexes that are not polyhedra, and discuss them first. We follow \cite[Ch. 7E]{arp}.

For a regular polyhedron $\K$ in $\mathbb{E}^3$ with symmetry group $G(\K)$, each subgroup $G_j$ of $G(\K)$ has order $2$ and is generated by a reflection $R_j$ in a point, line, or plane (a reflection in a line is a half-turn about the line). Thus $G(\K)$ is generated by $R_0,R_1,R_2$. We let $dim(R_j)$ denote the dimension of the mirror (fixed point set) of the reflection $R_j$ for each $j$, and call the vector $(dim(R_0), dim(R_1),dim(R_2))$ the {\em complete mirror vector} of $\K$;  this is just the dimension vector of \cite[Ch. 7E]{arp}. The use of the qualification ``complete" will become clear in the next section. The {\em distinguished generators\/} $R_0,R_1,R_2$ of $G(\K)$ satisfy (at least) the Coxeter-type relations
\begin{equation}
\label{relone}
R_{0}^{2} = R_{1}^{2} = R_{2}^{2} =
(R_{0}R_{1})^{p} = (R_{1}R_{2})^{q} = (R_{0}R_{2})^{2} = I,
\end{equation}
the identity mapping, where $p$ and $q$ determine the {\em type} $\{p,q\}$ of $\K$.

The complete enumeration of the regular polyhedra naturally splits into four steps of varying degrees of difficulty:\ the finite polyhedra, the planar apeirohedra, the blended apeirohedra, and the pure (non-blended) apeirohedra. An {\it apeirohedron\/} is simply an infinite polyhedron.

There are just $18$ finite regular polyhedra: the five (convex) Platonic solids
\[ \{3,3\},\, \{3,4\},\, \{4,3\},\, \{3,5\},\, \{5,3\};\]
the four Kepler-Poinsot star-polyhedra
\[ {\textstyle\{3,\frac{5}{2}\},\, \{\frac{5}{2},3\},\, \{5,\frac{5}{2}\},\, \{\frac{5}{2},5\}}, \]
where faces and vertex-figures are planar, but are allowed to be star polygons;
and the Petrie-duals of these nine polyhedra. (Recall that the {\em Petrie dual\/} of a regular polyhedron $\mathcal{P}$ has the same vertices and edges as $\mathcal{P}$; however, its faces are the {\em Petrie polygons\/} of $\mathcal{P}$, whose defining property is that two successive edges, but not three, are edges of a face of $\mathcal{P}$. Thus the new faces are ``zig-zags", leaving a face of $\mathcal{P}$ after traversing two of its edges.)

The $6$ planar regular apeirohedra comprise the three familiar regular plane tessellations by squares, triangles, or hexagons,
\[ \{4,4\},\, \{3,6\},\, \{6,3\}, \]
and their Petrie-duals.

The remaining regular apeirohedra are genuinely $3$-dimensional and fall into two families.

There are exactly $12$ regular apeirohedra that in some sense are reducible and have components that are regular figures of dimensions $1$ and $2$ . These apeirohedra are {\it blends} of a planar regular apeirohedron, and a line segment $\{\,\}$ or linear apeirogon $\{\infty\}$. This explains why there are $12=6\cdot 2$ blended (or non-pure) aperiohedra. For example, the blend of the standard square tessellation $\{4,4\}$ and the infinite apeirogon $\{\infty\}$, denoted $\{4,4\}\#\{\infty\}$, is an apeirohedron whose faces are helical apeirogons (over squares), rising above the squares of $\{4,4\}$, such that $4$ meet at each vertex; the orthogonal projections of $\{4,4\}\#\{\infty\}$ onto their component subspaces recover the original components, the square tessellation and linear apeirogon.

Note that each blended polyhedron really represents an entire family of polyhedra of the same kind, where the polyhedra in a family are determined by a parameter describing the relative scale of the two component figures. Thus there are infinitely many polyhedra of each kind, up to similarity, and our original count really refers to the 12 kinds rather than individual polyhedra.

Finally there are $12$ regular apeirohedra that are irreducible, or {\it pure\/} (non-blended). In a sense, they fall into a single family, derived from the standard regular cubical tessellation. The $12$ polyhedra in this family naturally are interrelated by a net of geometric operations (on polyhedra) and algebraic operations (on symmetry groups), which include the following:\ the duality operation; the previously mentioned Petrie-operation (of passing to the Petrie-dual); the {\it facetting\/} operation (of replacing the faces of a regular polyhedron by its {\it holes\/}, which are edge paths that successively take the second exit on the right at each vertex, while keeping all the vertices and edges unchanged); two lesser known operations called {\it halving\/} and {\it skewing\/}; and certain combinations of these operations.

We list these $12$ pure apeirohedra in the following table taken from \cite[p. 225]{arp}, which also highlights the fact that there are just $12$ polyhedra of this kind.
\medskip

\begin{table}[htb]
\begin{center}
\begin{tabular}{|c||ccc|c|c|}
\hline
mirror vector& $\{3,3\}$ & $\{3,4\}$ & $\{4,3\}$  &faces&vertex-fig.\\
\hline  \hline
(2,1,2) \,&\, $\{6,6 | 3\}$ \,&\, $\{6,4 | 4\}$ \,&\, $\{4,6 | 4\}$ \,&\,planar\,&\,skew\,  \\[.05in]
(1,1,2) \,&\, $\{\infty,6\}_{4,4}$ \,&\, $\{\infty,4\}_{6,4}$ \,&\, $\{\infty,6\}_{6,3}$ \,&\,helical\,&\,skew\, \\[.05in]
(1,2,1) \,&\, $\{6,6\}_{4}$ \,&\, $\{6,4\}_{6}$ \,&\, $\{4,6\}_{6}$\,&\,skew\,&\,planar\, \\[.05in]
(1,1,1) \,&\, $\{\infty,3\}^{(a)}$ \,&\, $\{\infty,4\}_{\cdot,*3}$ \,&\,
$\{\infty,3\}^{(b)}$ \,&\,helical\,&\,planar\, \\[.05in]
\hline
\end{tabular}
\end{center}
\caption{The 12 pure apeirohedra in $\mathbb{E}^3$}\label{tabzero}
\end{table}
\medskip

\noindent
In this table, the first column gives the complete mirror vector, and the last two describe if the faces and vertex-figures are planar, skew or helical regular polygons; the geometric nature of the faces and vertex-figures only depends on the mirror vector. The second, third, and fourth column are indexed by the finite Platonic polyhedra whose rotation or full symmetry group is intimately related to the special group.

The three polyhedra along the top row are the famous Petrie-Coxeter polyhedra, which along with those in the third row comprise the pure regular polyhedra with finite faces. The pure polyhedra with infinite, helical faces are listed in the second and last row; those in the last row occur in two enantiomorphic (mirror image) forms, since their symmetry group is generated by half-turns and consists only of proper isometries. The fine Schl\"afli symbols for the polyhedra in the table signify defining relations for the symmetry groups; for example, extra relations often specify the orders of the elements $R_{0}R_{1}R_{2}$, $R_{0}R_{1}R_{2}R_{1}$ or $R_{0}(R_{1}R_{2})^{2}$. These orders correspond to the lengths of the Petrie paths, of the holes (paths traversing edges where the new edge is chosen to be the second on the right according to some local orientation), and of the 2-zigzags (paths traversing edges where the new edge is chosen to be the second on the right, but reversing orientation on each step).

\begin{figure}[htb]
\begin{center}
\includegraphics[width=1.8in]{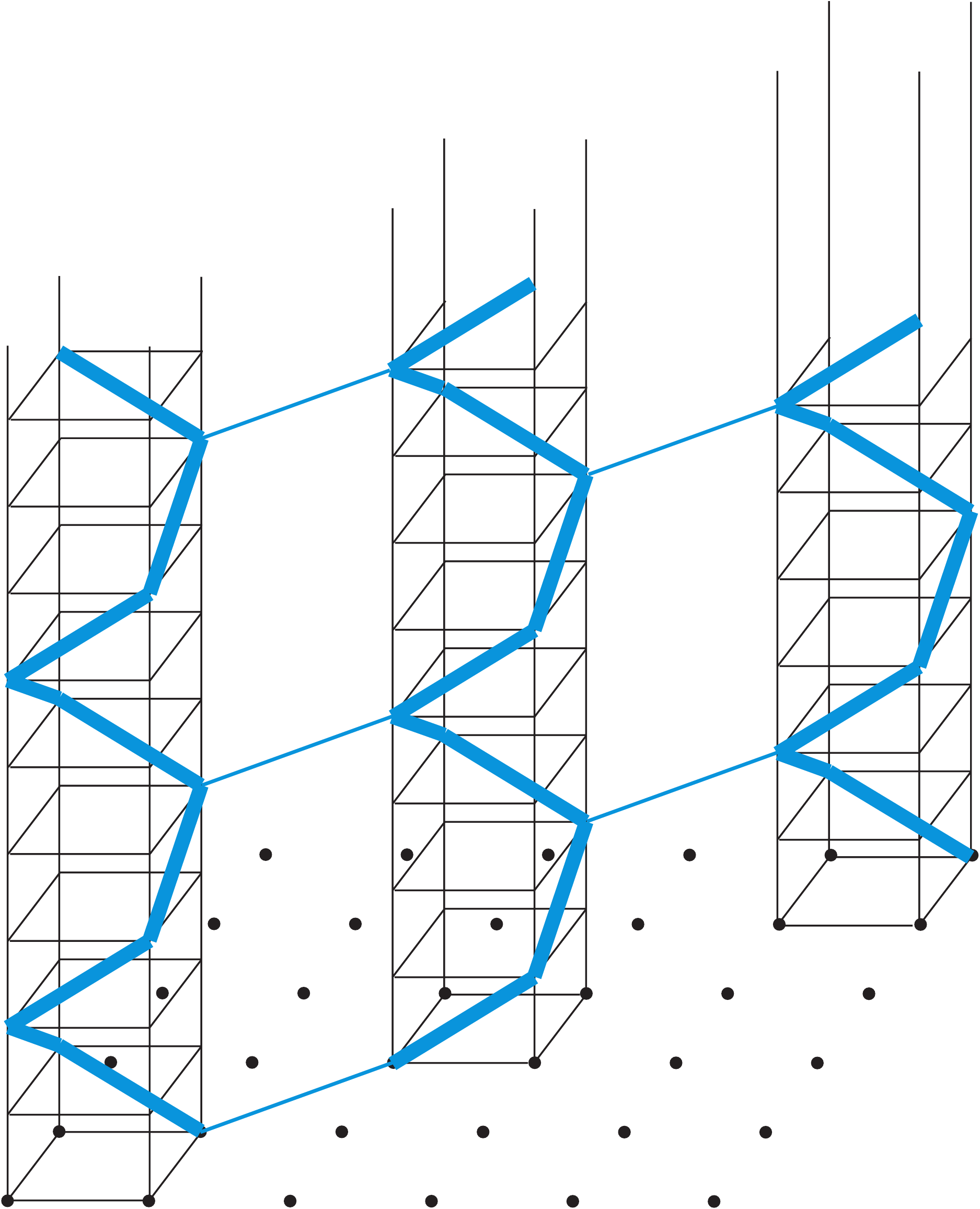}\qquad\qquad
\includegraphics[width=1.8in]{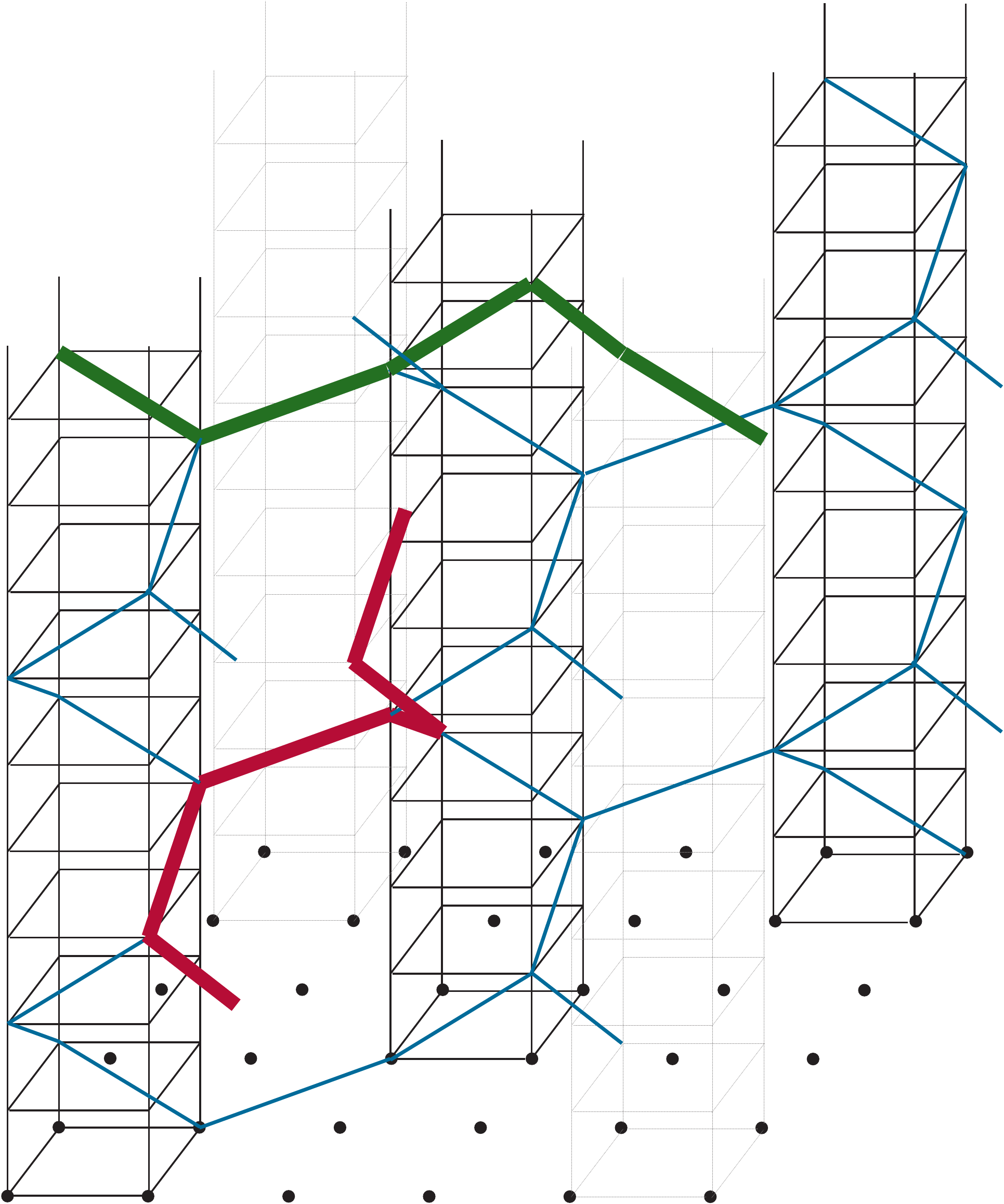}\\[.2in]
\end{center}
\caption{The helix-faced regular polyhedron $\{\infty,3\}^{(b)}$, with symmetry group requiring the single extra relation $(R_{0}R_{1})^{4}(R_{0}R_{1}R_{2})^{3} = (R_{0}R_{1}R_{2})^{3} (R_{0}R_{1})^{4}$.}
\label{figone}
\end{figure}

The regular polyhedron $\{\infty,3\}^{(b)}$ is illustrated in Figure~\ref{figone}; three helical faces meet at each vertex. Some faces have a vertical axis; they are helices over squares, like the ones shown on the left, and are joined by horizontal edges. The remaining faces have axes parallel to the remaining two coordinate axes; one copy of each is shown on the right.

In summary we have

\begin{theorem}
\label{polyhs}
There are precisely $48$ regular polyhedra in $\mathbb{E}^3$, up to similarity and scaling of components (when applicable). The  list comprises $18$ finite polyhedra and $30$ apeirohedra.
\end{theorem}

\section{Non-simply flag-transitive complexes}
\label{nonsim}

In order to complete the classification of regular polygonal complexes in $\mathbb{E}^3$ it remains to consider  complexes with three or more faces around each edge. For convenience we split the discussion into two cases according to the size of the flag stabilizers. Throughout this and the next section we follow \cite{pelsch1}, \cite{pelsch2}.

Quite surprisingly, up to similarity, there are just four regular polygonal complexes that are not simply flag-transitive. They can be characterized as the regular complexes $\K$ that occur as $2$-skeletons of regular $4$-apeirotopes $\mathcal{P}$ in $\mathbb{E}^3$ (see~\cite[Ch. 7F]{arp}). The {\em $2$-skeleton} of a $4$-apeirotope is the incidence structure determined by its vertices, edges and polygons. These $4$-apeirotopes in $\mathbb{E}^3$ are, by definition, the discrete faithful realizations of abstract regular polytopes of rank $4$ in $\mathbb{E}^3$, so their combinatorial rank is $1$ higher than the dimension of the ambient space.

There are precisely eight regular $4$-apeirotopes $\mathcal{P}$ in $\mathbb{E}^3$, occurring in pairs of Petrie-duals as shown in (\ref{4apeirotopes}). The Petrie-dual of a regular $4$-apeirotope $\mathcal{P}$ is obtained by replacing the distinguished involutory generators $T_0,T_1,T_2,T_3$ of its symmetry group $G(\mathcal{P})$ by the new involutory generators
\[ T_0,T_1T_3,T_2,T_3 \]
of $G(\mathcal{P})$, and then applying Wythoff's construction with these new generators and with the same initial vertex as for $\mathcal{P}$ itself. Every pair of Petrie-duals contributes just one regular polygonal complex $\K$, since Petrie-duals have isomorphic $2$-skeletons. Thus there are just four such complexes $\K$. Two of the eight apeirotopes $\mathcal{P}$ have (finite convex) square $2$-faces, $4$ occurring at each edge; and six have (infinite planar) zigzag $2$-faces, with either $3$ or $4$ at each edge. Our notation follows \cite[Ch. 7F]{arp}.
\begin{equation}
\label{4apeirotopes}
\begin{array}{ccc}
\{4, 3, 4\} &\quad &\{\{4, 6 \,|\,4\}, \{6, 4\}_3\}     \\[.04in]
\{\{\infty, 3\}_6 \# \{ \, \}, \{3, 3\}\}
&\quad &\{\{\infty, 4\}_4 \# \{\infty\}, \{4, 3\}_3\}  \\[.04in]
\{\{\infty, 3\}_6 \# \{ \, \}, \{3, 4\}\}
&\quad &\{\{\infty, 6\}_3 \# \{\infty\}, \{6, 4\}_3\}  \\[.04in]
\{\{\infty, 4\}_4 \# \{\, \}, \{4, 3\}\}
&\quad &\{\{\infty, 6\}_3 \# \{\infty\}, \{6, 3\}_4\}  \\[.04in]
\end{array}
\end{equation}
The two apeirotopes in the top row are the standard cubical tessellation $\{4,3,4\}$ in $\mathbb{E}^3$; and its Petrie-dual $\{\{4,6\,|\,4\},\{6,4\}_3\}$, whose rank 3 faces are Petrie-Coxeter polyhedra $\{4,6\,|\,4\}$ and whose vertex-figures are Petrie-duals $\{6, 4\}_3$ of octahedra $\{3,4\}_6$. The $2$-skeleton of the cubical tessellation is the simplest regular polygonal complex that is not simply flag-transitive.

The other six apeirotopes have finite crystallographic regular polyhedra as vertex-figures, namely either tetrahedra $\{3,3\}$, octahedra $\{3,4\}$, or cubes $\{4,3\}$, or Petrie-duals of one of those; their rank 3 faces are blends, namely of the Petrie-duals $\{\infty,3\}_6$ or $\{\infty,4\}_4$ of the plane tessellations $\{6,3\}$ or $\{4,4\}$, respectively, with the line segment $\{ \,\}$ or linear apeirogon $\{\infty\}$ (see \cite[Ch.7E]{arp}).

The number of faces $r$ around an edge of the $2$-skeleton $\K$ is just the last entry in the Schl\"afli symbol (the basic symbol $\{p,q,r\}$) of the underlying $4$-apeirotope $\mathcal{P}$ (or, equivalently, of the Petrie dual of $\mathcal{P}$). Hence, $r=4$, $3$, $4$ or $3$, respectively.

Among the regular polygonal complexes $\K$, the non-simply transitive complexes can also be characterized as those that have face mirrors. A {\em face mirror\/} of $\K$ is an affine plane in $\mathbb{E}^3$ that contains a face of $\K$ and is the mirror of a plane reflection in $G(\K)$. Clearly, regular complexes $\K$ with face mirrors must have planar faces, and every face must span a face mirror; moreover, the plane reflection in a face mirror of $\K$ fixes every flag of $\K$ that lies in this face mirror, and hence generates the corresponding flag stabilizer.

In summary we have

\begin{theorem}
\label{Petskels}
Up to similarity, there are just four non-simply flag-transitive regular polygonal complexes in $\mathbb{E}^3$, each given by the common $2$-skeleton of the two regular $4$-apeirotopes from a pair of Petrie-duals. These infinite complexes are precisely the regular polygonal complexes in $\mathbb{E}^3$ that have face mirrors.
\end{theorem}

\section{Simply flag-transitive complexes}
\label{simco}

The class of simply flag-transitive regular polygonal complexes in $\mathbb{E}^3$ is much richer and comprises all finite or infinite regular polyhedra. As we have already described the regular polyhedra in Section \ref{regpol}, we can confine ourselves here to those complexes that are not polyhedra.

Thus let $\K$ be an (infinite) simply flag-transitive complex such that $G=G(\K)$ is affinely irreducible, let $r\geq 3$, and let $\{F_0,F_1,F_2\}$ denote the base flag. Then the two subgroups $G_0$ and $G_1$ of $G$ are again of order $2$, and are generated by some point, line, or plane reflection $R_0$, and some line or plane reflection $R_1$, respectively; however, $G_2$ is a cyclic or dihedral group of order $r$. The {\em mirror vector} $(dim(R_0), dim(R_1))$ of $\K$ now has only two components recording the dimensions $dim(R_0)$ and $dim(R_1)$ of the mirrors of $R_0$ and $R_1$, respectively. (For polyhedra, $G_2$ is also generated by a reflection, and the complete mirror vector records the dimensions of all three mirrors.)

The vertex-stabilizer subgroup $G_{F_0}$ in $G$ of the base vertex $F_0$ is called the {\em vertex-figure group\/} of $\K$ at $F_0$, and is a finite group since $\K$ is discrete. In particular, $G_{F_0}=\langle R_1, G_2 \rangle$, and $G_{F_0}$ acts simply flag-transitively on the graph that forms the vertex-figure of $\K$ at $F_0$. (A flag of a graph is just an incident vertex-edge pair.)  Similarly, the face-stabilizer $G_{F_2}$ in $G$ of the base face $F_2$ is given by $G_{F_2}=\langle R_0, R_1 \rangle$ and is isomorphic to a (finite or infinite) dihedral group acting simply transitively on the flags of $\K$ containing $F_2$.

The enumeration of the simply flag-transitive regular complexes for a given mirror vector is typically rather involved. A good number of complexes must be discovered by direct geometric or algebraic methods. Others then can be derived by operations applied to these complexes; that is, the new complexes are obtained by suitably modifying $R_0$ and $R_1$ while keeping the base vertex and preserving the group $\langle G_0, G_1, G_2 \rangle$ as a (possibly proper) subgroup of symmetries. In this vein, the explicit enumeration of the simply flag-transitive complexes begins in \cite{pelsch1} with the determination of the complexes with mirror vector $(1,2)$, and then proceeds in \cite{pelsch2} with the description of those for the remaining mirror vectors, accomplished by a mix of direct methods, applications of operations, and elimination of certain cases. At the end, we arrive at the following theorem.

\begin{theorem}
\label{fullclassifsim}
Up to similarity, there are exactly $21$ simply flag-transitive regular polygonal complexes in $\mathbb{E}^3$ that are not regular polyhedra.
\end{theorem}

Thus, counting also the regular polyhedra from Theorem~\ref{polyhs}, there is total of $69$ simply flag-transitive regular complexes, up to similarity and scaling of components for blended polyhedra.

Table~\ref{tabone} lists the $21$ simply flag-transitive complexes by mirror vector, and records their data concerning the pointwise edge stabilizer $G_2$, the number $r$ of faces surrounding an edge, the structure of the faces and vertex-figures, the vertex-set, and the structure of the special group $G_*$. In the face column we have used the symbols $p_c$, $p_s$, $\infty_2$, or $\infty_k$ with $k=3$ or $4$, respectively, to indicate that the faces are {\em \underbar{c}onvex} $p$-gons, {\em \underbar{s}kew} $p$-gons, planar zigzags, or helical polygons over {\em $k$-gons}. (A planar zigzag is viewed as a helix over a $2$-gon, hence our notation. Clearly, the subscript in $3_c$ is redundant.)  We also set
\[ V_{a}:=a\mathbb{Z}^{3}\!\setminus\! ((0,0,a)\!+\!\Lambda_{(a,a,a)}),\;\;
W_{a}:= 2\Lambda_{(a,a,0)} \cup ((a,-a,a)\!+\!2\Lambda_{(a,a,0)}), \]
to have a short symbol available for the vertex-sets of some complexes. The vertex-figures of polygonal complexes are finite geometric graphs, so an entry in the vertex-figure column describing a solid figure is meant to represent the edge-graph of this figure, with ``double" indicating the double edge-graph. The abbreviation ``ns-cuboctahedron" stands for the edge graph of a certain ``non-standard cuboctahedron", a realization in $\mathbb{E}^3$ of the (abstract) cuboctahedron with non-planar square faces.

\begin{table}[htb]
\begin{center}
{\begin{tabular}{|c|c|c|c|c|c|c|c|}
\hline
$\,$mirror vector$\,$ & $\;$complex$\;$&$\;G_2\;$ & $\;r\;$ &$\;$face $\;$&$\;$vertex-figure$\;$ &$\;$vertex-set$\;$ & $\;$special group$\;$\\[.05in]
\hline
\hline
$(1,2)$  & $\K_1(1,2)$& $D_2$  & $4$ &$4_s$ & cuboctahedron&$\Lambda_{(a,a,0)}$&$[3,4]$\\
\hline
             & $\K_2(1,2)$& $C_3$& $3$ &$4_s$ & cube&$\Lambda_{(a,a,a)}$&$[3,4]$\\
\hline
             & $\K_3(1,2)$& $D_3$& $6$   &$4_s$ &double cube&$\Lambda_{(a,a,a)}$&$[3,4]$\\
\hline
             & $\K_4(1,2)$& $D_2$& $4$ & $6_s$&octahedron&$a\mathbb{Z}^3$&$[3,4]$\\
\hline
             & $\K_5(1,2)$& $D_2$& $4$ &$6_s$&double square&$V_a$&$[3,4]$\\
\hline
             & $\K_6(1,2)$& $D_4$& $8$ &$6_s$&double octahedron&$a\mathbb{Z}^3$&$[3,4]$\\
\hline
             & $\K_7(1,2)$& $D_3$& $6$ &$6_s$&double tetrahedron&$W_a$&$[3,4]$\\
\hline
             & $\K_8(1,2)$& $D_2$& $4$ &$6_s$&cuboctahedron&$\Lambda_{(a,a,0)}$&$[3,4]$\\
\hline
\hline
$(1,1)$  & $\K_1(1,1)$& $D_3$& $6$ &$\infty_3$&double cube&$\Lambda_{(a,a,a)}$&$[3,4]$\\
\hline
             & $\K_2(1,1)$& $D_2$& $4$ &$\infty_3$&double square&$V_a$&$[3,4]$\\
\hline
             & $\K_3(1,1)$& $D_4$& $8$ &$\infty_3$&double octahedron&$a\mathbb{Z}^3$ & $[3,4]$ \\
\hline
             & $\K_4(1,1)$& $D_3$ & $6$ &$\infty_4$& double tetrahedron&$W_a$ &$[3,4]$\\
\hline
             & $\K_5(1,1)$& $D_2$& $4$ &$\infty_4$&ns-cuboctahedron&$\Lambda_{(a,a,0)}$&$[3,4]$ \\
\hline
             & $\K_6(1,1)$& $C_3$& $3$ &$\infty_4$&tetrahedron& $W_a$&$[3,4]^+$\\
\hline
             & $\K_7(1,1)$& $C_4$& $4$ &$\infty_3$&octahedron&$a\mathbb{Z}^3$ & $[3,4]^+$ \ \\
\hline
             & $\K_8(1,1)$& $D_2$& $4$ &$\infty_3$&ns-cuboctahedron&$\Lambda_{(a,a,0)}$&$[3,4]$ \\
\hline
             & $\K_9(1,1)$& $C_3$& $3$ &$\infty_3$&cube&$\Lambda_{(a,a,a)}$&$[3,4]^+$ \\
\hline
\hline
$(0,1)$  & $\K(0,1)$& $D_2$& $4$ &$\infty_2$&ns-cuboctahedron&$\Lambda_{(a,a,0)}$&$[3,4]$ \\
\hline
\hline
$(0,2)$  & $\K(0,2)$& $D_2$& $4$ &$\infty_2$&cuboctahedron& $\Lambda_{(a,a,0)}$&$[3,4]$ \\
\hline
\hline
$(2,1)$  & $\K(2,1)$& $D_2$& $4$ & $6_c$ &ns-cuboctahedron&$\Lambda_{(a,a,0)}$&$[3,4]$ \\
\hline
\hline
$(2,2)$  & $\K(2,2)$& $D_2$& $4$ &$3_c$&cuboctahedron&$\Lambda_{(a,a,0)}$&$[3,4]$\\
\hline
\end{tabular}}
\end{center}
\medskip
\caption{The 21 simply flag-transitive regular complexes in $\mathbb{E}^3$ which are not regular polyhedra.}
\label{tabone}
\end{table}

As an example, the faces of the complex $\K_6(1,2)$ are the Petrie polygons of all cubes of the cubical tessellation of $\mathbb{E}^3$; so in particular, the vertices and edges of $\K_6(1,2)$, respectively, comprise all vertices and edges of the cubical tessellation. Recall that every edge of a cube belongs to precisely two Petrie polygons of the same cube. Since every edge belongs to four cubes in the cubical tessellation, every edge must belong to eight Petrie polygons of cubes in $\K_6(1,2)$. The complex $\K_4(1,2)$ is a proper subcomplex of $\K_6(1,2)$ obtained by taking only the Petrie polygons of alternate cubes.
The complex $\K_5(1,2)$ is another subcomplex of $\K_6(1,2)$ consisting only of the Petrie polygons with vertices in the set $V_a$ defined above.

\section{Chiral polyhedra}
\label{chirpoly}

Chiral polyhedra in $\mathbb{E}^3$ are the most interesting kind of nearly regular polyhedra; their geometric symmetry groups have two orbits on the flags, such that adjacent flags are in distinct orbits.

The structure results for the symmetry groups of regular polygonal complexes carry over to chiral polyhedra as follows (see \cite{chiral1,chiral2}). Let $\K$ be a chiral polyhedron in $\mathbb{E}^3$ with symmetry group $G=G(\K)$, let $\Phi := \{F_{0},F_{1},F_{2}\}$ be a base flag of $\K$, and let $F_{0}',F_{1}',F_{2}'$ denote the faces of $\K$ with $F'_0<F_1$, $F_0 < F'_1 < F_2$, $F_1<F'_2$  and $F'_{j} \neq F_j$ for $j=0,1,2$. Then $G$ is generated by symmetries $S_{1},S_{2}$  of $\K$, called the {\em distinguished generators\/} of $G$ (relative to $\Phi$), where $S_{1}$ leaves the base face $F_2$ invariant and cyclically permutes the vertices of $F_2$ such that $F_{1}S_{1} = F'_{1}$ (and thus $F'_{0}S_{1} = F_{0}$), and $S_{2}$ leaves the base vertex $F_0$ invariant and cyclically permutes the vertices in the vertex-figure at $F_0$ such that $F_{2}S_{2} = F'_{2}$ (and thus $F'_{1}S_{2} = F_{1}$). Then, in analogy to~(\ref{relone}),
\begin{equation}
\label{reltwo}
S_{1}^p = S_{2}^q = (S_{1}S_{2})^{2}  = I,
\end{equation}
where $\{p,q\}$ is the Schl\"afli type of $\K$. The involutory symmetry $T := S_{1}S_{2}$ interchanges the two end vertices of $F_{1}$ as well as the two faces meeting at $F_{1}$; that is, combinatorially, $T$ acts like a half-turn about the midpoint of an edge. This symmetry $T$ plays a critical role, in that it allows to employ a variant of Wythoff's construction (see \cite{coxeter}) to  reconstruct a chiral polyhedron from its symmetry group.

Note that the symmetry groups of regular polyhedra in $\mathbb{E}^3$ have a subgroup of index at most $2$ with properties very similar to those of the group of a chiral polyhedron. In fact, if $\mathcal{P}$ is a regular polyhedron and $R_{0},R_{1},R_{2}$ are the distinguished generators of its symmetry group $G(P)$ (relative to $\Phi$), then $\widehat{S}_{1} := R_{0}R_{1}$ and $\widehat{S}_{2} := R_{1}R_{2}$ generate the {\em combinatorial rotation subgroup\/}, or {\em even subgroup\/}, $G^+(\mathcal{P}) := \langle\widehat{S}_{1},\widehat{S}_{2}\rangle$ of $G(P)$, of index $1$ or $2$. Now $\widehat{T} := \widehat{S}_{1}\widehat{S}_{2} = R_{0}R_{2}$ has properties similar to $T$. Whenever $G^+(\mathcal{P})$ has index $2$ in $G(\mathcal{P})$ we say that $\mathcal{P}$ is {\em directly regular} or {\em orientable}.

Combinatorially speaking, chiral polyhedra have maximal ``rotational" symmetry but no ``reflexive" symmetry. (This does not mean that $S_1$ and $S_2$ are actually geometric rotations! )  Thus our term ``chiral" really means ``maximal chiral". By contrast, again combinatorially speaking, regular polytopes have maximal ``reflexive" symmetry. (Here $R_0,R_1,R_2$ are actually reflections, in points, lines, or planes.)

Chirality, in this sense of ``maximal chirality", does not make any appearance in the classical theory of highly-symmetric figures in Euclidean spaces. This may explain why chiral polyhedra were only described and enumerated quite recently, in \cite{chiral1,chiral2}.

The complete classification starts off with the observation that chiral polyhedra are necessarily pure apeirohedra; that is, infinite polyhedra that are not naturally ``blends" of two lower-dimensional structures, and hence have an affinely irreducible symmetry group. In short, unlike regular polyhedra, chiral polyhedra can neither be finite nor planar or blended.

\begin{table}[htb]
\begin{center}
\begin{tabular}{|c|c|c|c|}
\hline
Type &$\{6,6\}$ &$\{4,6\}$  &$\{6,4\}$\\
\hline
\hline
&&&\\[-.1in]
Notation & $P(a,b)$& $Q(c,d)$& $Q(c,d)^*$\\[.02in]
\hline
&&&\\[-.1in]
Parameters&$a,b\in\mathbb{Z},\,(a,b)=1$&$c,d\in\mathbb{Z},\,(c,d)=1$&$c,d\in\mathbb{Z},\,(c,d)=1$ \\[.02in]
\hline
&&&\\[-.1in]
Chiral &$b\neq \pm a$&$c,d\neq 0$&$c,d\neq 0$\\[.04in]
\hline
&&&\\[-.1in]
Regular&$P(a,\!-a)\!=\!\{6,\!6\}_{4}$  &$Q(a,\!0)\!=\!\{4,\!6\}_{6}$&$Q(a,\!0)^*\!=\!\{6,\!4\}_{6}$  \\[.02in]
Polyhedra &$\quad\,P(a,\!a)\!=\!\{6,\!6|3\}$&$\;\,Q(0,\!a)\!=\!\{4,\!6|4\}$&$\;\,Q(0,\!a)^*\!=\!\{6,\!4|4\}$  \\[.03in]
\hline
&geom.$\!$ self-dual,&&\\[.01in]
&$P(a,b)^*\cong P(a,b)$&&\\[.02in]
\hline
&&&\\[-.1in]
Special Group& $[3,3]^{+}\times \langle -I\rangle$& $[3,4]$&$[3,4]$\\[.03in]
\hline
\end{tabular}
\end{center}
\medskip
\caption{The finite-faced chiral polyhedra, along with their related regular polyhedra. }
\label{tabtwo}
\end{table}

The classification of chiral apeirohedra is quite elaborate and naturally breaks down into analyzing the finite-faced and the helix-faced polyhedra (see \cite{chiral1,chiral2}). The possible apeirohedra fall into six infinite $2$-parameter families (up to congruence).  In each family, all but two polyhedra are chiral; the two exceptional polyhedra are regular and are among those described in Section~\ref{regpol}. Tables~\ref{tabtwo} and~\ref{tabthree} list the families of polyhedra by Schl\"afli type, along with the two regular polyhedra occurring in each family; in the three families in Table~\ref{tabthree}, one exceptional polyhedron is finite. Also included is data about the {\em special group\/} of a polyhedron, that is, the quotient of the geometric symmetry group by its translation subgroup; here $[3,3]^+$ and $[3,4]^+$ denote the tetrahedral or octahedral rotation group, respectively, and $[3,4]$ the full octahedral group.

It is quite remarkable that a regular polyhedron cannot have both skew faces and skew vertex-figures. However, finite-faced chiral polyhedra must necessarily have both skew faces and skew vertex-figures. In fact, the generators $S_1,S_2$ of the symmetry group must be rotatory reflections in this case, resulting in skew faces and skew vertex-figures. Note, however, that the rotation subgroups for the regular polyhedra occurring in the three families of finite-faced polyhedra of Table~\ref{tabtwo} also have
generators $S_1,S_2$ which are rotatory reflections, but here the position of the base vertex forces planarity of faces or vertex-figures.

\begin{figure}[htb]
\begin{center}
\ \includegraphics[width=1.6in]{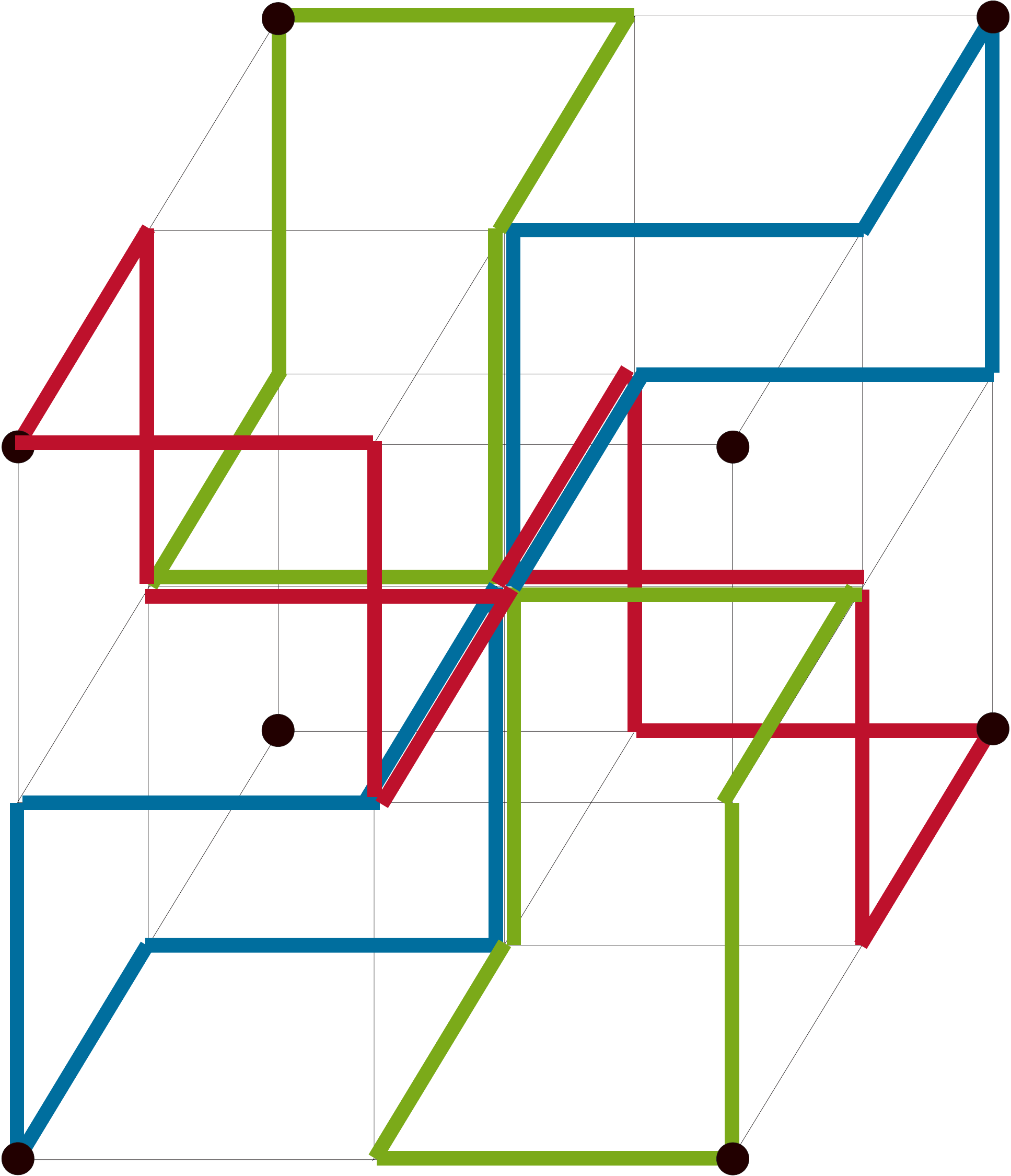}\qquad
\includegraphics[width=2in]{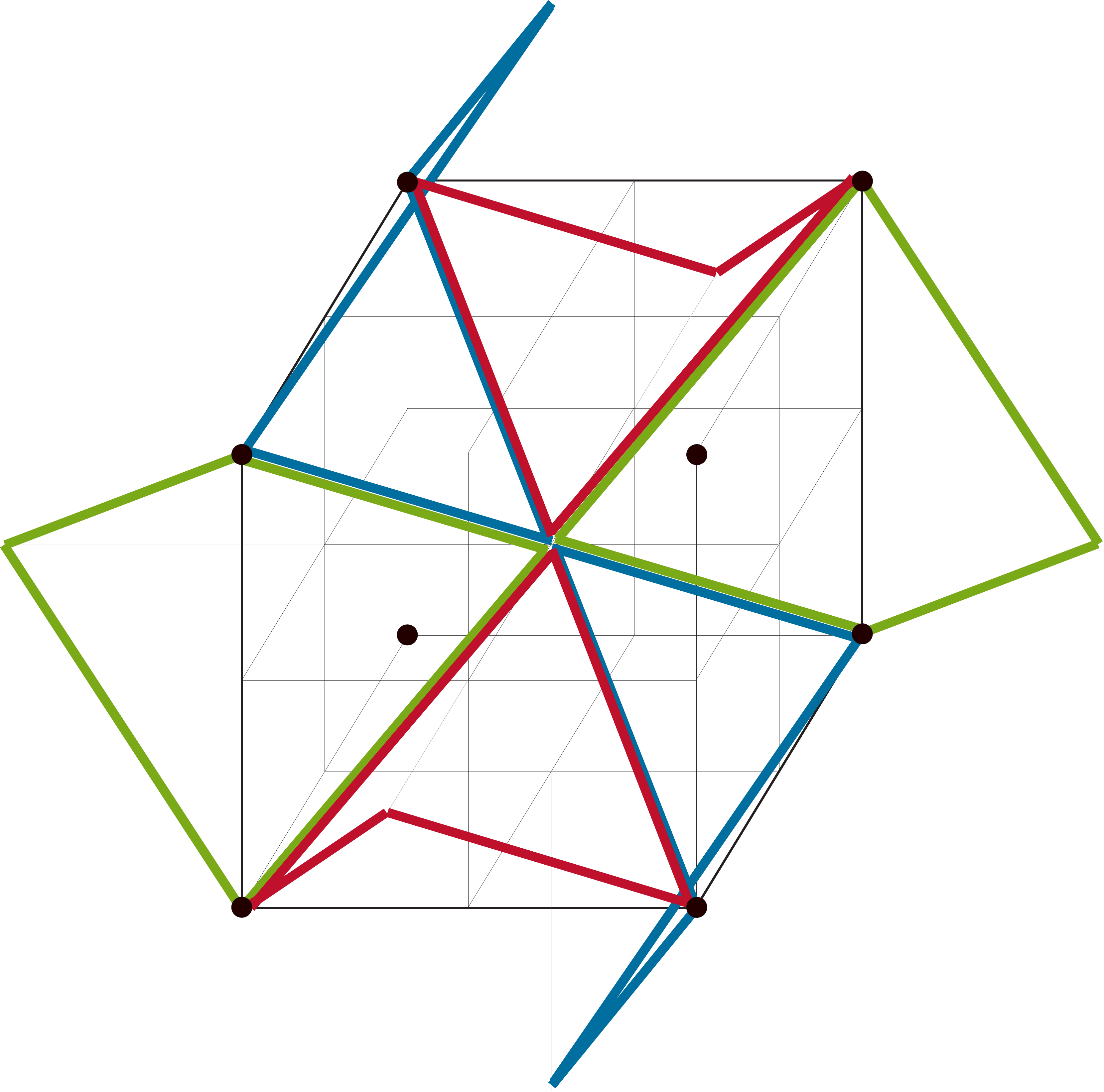}\\[.2in]
\end{center}
\caption{The finite-faced chiral apeirohedra $P(1,0)$ and $Q(1,1)$, of types $\{6,6\}$ and $\{4,6\}$, respectively. Depicted is the neighborhood of a single vertex, where $6$ skew hexagonal faces or $6$ skew square faces meet. Each apeirohedron expands in a consistent manner throughout space such that all vertex neighborhoods are congruent to the one shown.}
\label{figtwo}
\end{figure}

Chiral apeirohedra with infinite faces must necessarily have helical faces spiraling over triangles or squares, as well as planar vertex-figures. The symmetry group is generated by a screw motion $S_1$ and a rotation $S_2$ in this case. Chiral helix-faced polyhedra unravel, in a sense, a ``crystallographic" Platonic polyhedron, namely the finite regular polyhedron in their respective family.

\begin{table}[htb]
\begin{center}
\begin{tabular}{|c|c|c|c|}
\hline
Type &$\{\infty,3\}$ &$\{\infty,3\}$&$\{\infty,4\}$\\
\hline
\hline
&&&\\[-.1in]
Notation &$P_1(a,b)$ & $P_2(c,d)$ & $P_3(c,d)$\\[.02in]
\hline
&&&\\[-.1in]
Parameters& $a,b\in\mathbb{R}$& $c,d\in\mathbb{R}$&$c,d\in\mathbb{R}$ \\[.04in]
\hline
&&&\\[-.1in]
Chiral & $b\neq \pm a$ & $c,d\neq 0$ & $c,d\neq 0$\\
\hline
&&&\\[-.1in]
Regular Polyhedra&$\;\;P_{1}(1,-1)=\{\infty,3\}^{(a)}$&$\quad\;P_{2}(1,0)=\{\infty,3\}^{(b)}$&$\quad\;\,P_{3}(0,1)=\{\infty,4\}_{\cdot, *3}$  \\[.03in]
&$P_{1}(1,1)=\{3,3\}$&$P_{2}(0,1)=\{4,3\}$& $P_{3}(1,0)=\{3,4\}$ \\[.03in]
\hline
&&&\\[-.1in]
Helices over& triangles& squares& triangles\\[.04in]
\hline
&&&\\[-.1in]
Special Group& $[3,3]^+$& $[3,4]^+$&$[3,4]^+$\\[.04in]
\hline
\end{tabular}
\end{center}
\caption{The helix-faced chiral polyhedra, along with their related regular polyhedra. }
\label{tabthree}
\end{table}

The regular polyhedra listed in Tables~\ref{tabtwo} and \ref{tabthree} comprise nine of the twelve pure regular apeirohedra in $\mathbb{E}^3$, namely those listed in Table~\ref{tabzero} with complete mirror vectors $(1,2,1)$, $(1,1,1)$ or $(2,1,2)$, as well as the three crystallographic Platonic polyhedra. The three remaining pure regular apeirohedra $\{\infty,6\}_{4,4}$, $\{\infty,4\}_{6,4}$ and $\{\infty,6\}_{6,3}$ all have complete mirror vector $(1,1,2)$ and do not occur in families alongside chiral polyhedra.

\begin{figure}[htb]
\begin{center}
\includegraphics[width=2in]{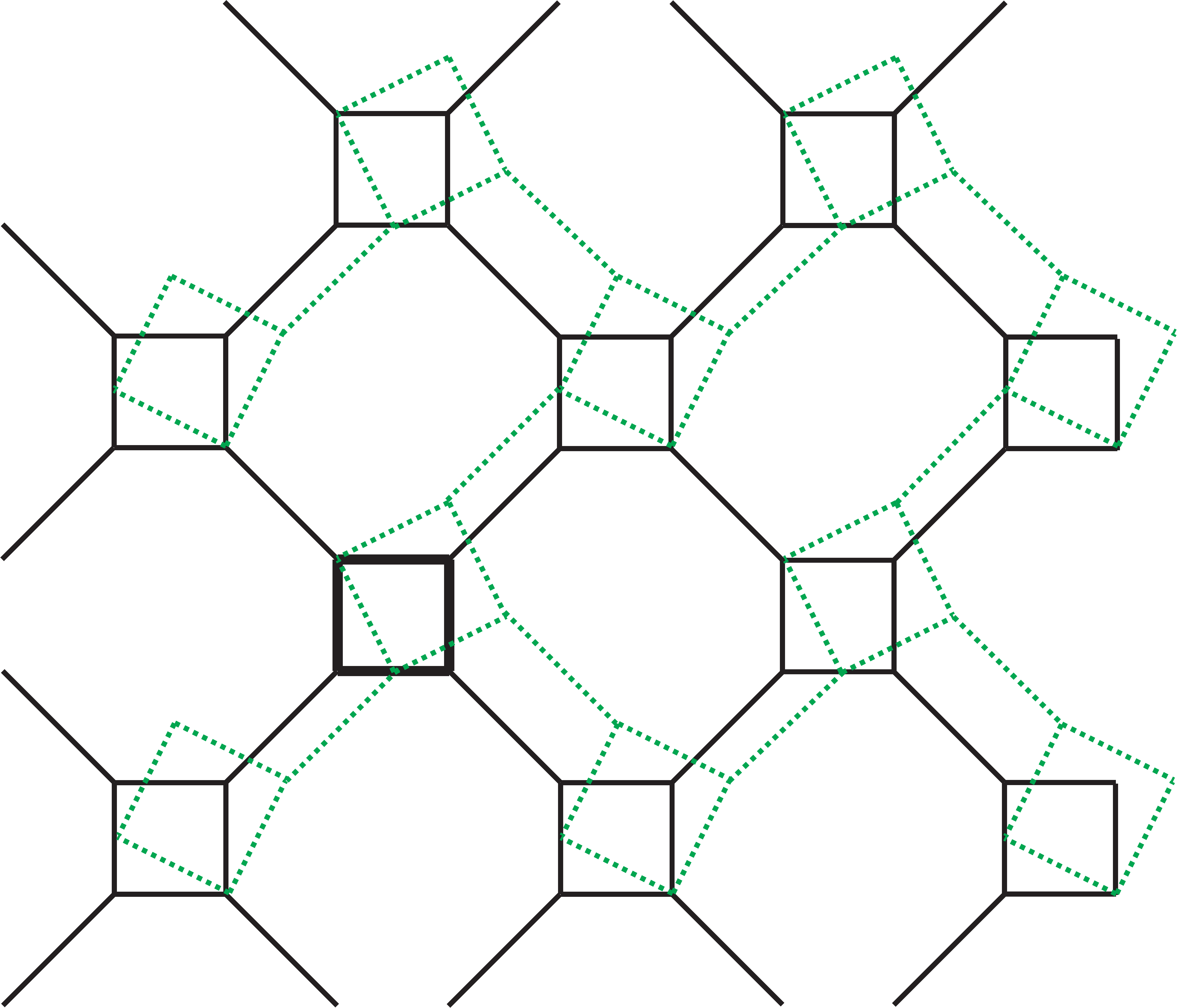}\qquad
\includegraphics[width=2in]{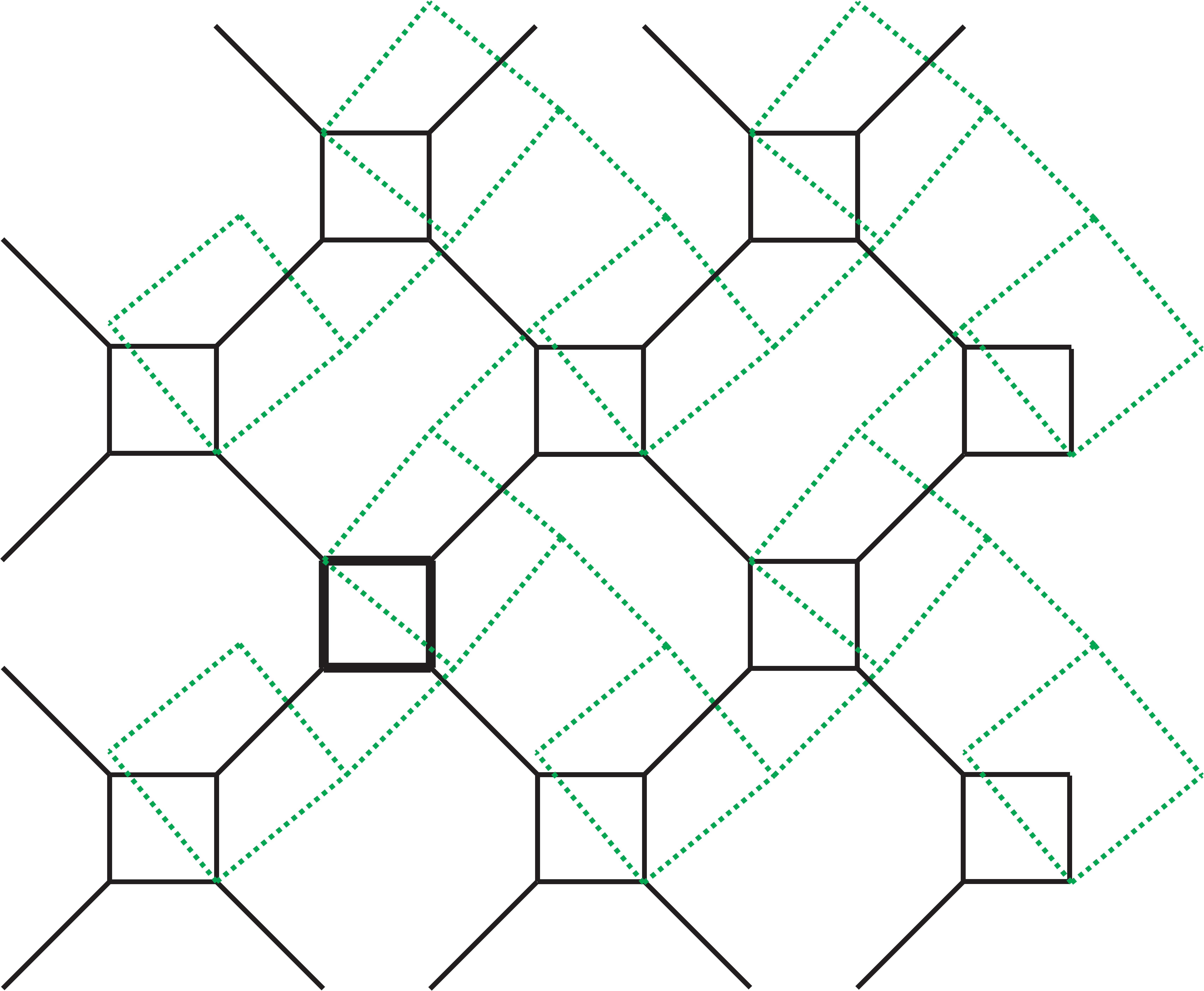}
\end{center}
\caption{The helix-faced polyhedron $P_{2}(1,0)$ and its deformations $P_{2}(1,d)$. The solid black or green dotted lines show the projection of the entire polyhedron $P_{2}(1,0)$ or $P_{2}(1,d)$, respectively, onto a ``horizontal" plane perpendicular to the axis of a ``vertical" helical face. The vertical helical faces of $P_{2}(1,0)$ or $P_{2}(1,d)$, respectively, then project onto the small black squares or small green squares; one black square, resulting from one vertical helical face of $P_{2}(1,0)$, is emphasized. As the parameter $d$ is changed continuously, the vertical and ``horizontal" helical faces move in such a way that the
axes of corresponding faces remain parallel throughout the process. Accordingly, the projections of $P_{2}(1,d)$ move continuously as well. The figures on the left and right, respectively, show projections of $P_{2}(1,d)$ when $d$ is small or when $d$ gets larger.}
\label{figcont}
\end{figure}

These six families of chiral (or regular) polyhedra have some amazing properties. For example, any two distinct finite-faced polyhedra of the same type are combinatorially non-isomorphic. In fact, $P(a,b)$ and $P(a',b')$ are isomorphic if and only if $(a',b')=\pm (a,b),\pm (b,a)$; and similarly, $Q(c,d)$ and $Q(c',d')$ are isomorphic if and only if $(c',d')=\pm (c,d),\pm (-c,d)$. Thus there are very many combinatorially distinct finite-faced chiral polyhedra. By contrast, as shown in Pellicer-Weiss~\cite{pelwei}, every helix-faced chiral polyhedron $P_1(a,b)$ or $P_2(c,d)$ is combinatorially isomorphic to the infinite regular polyhedron in its family. On the other hand, since the polyhedron $\{\infty, 4\}_{\cdot, *3}$ is not orientable, it cannot have chiral realizations. Every chiral polyhedron $P_3(c,d)$ is then isomorphic to the (combinatorial) orientable double cover of $\{\infty, 4\}_{\cdot, *3}$. Thus, up to isomorphism, there are just three helix-faced chiral polyhedra, each represented by a helix-faced regular polyhedron. But even more is true:\ in a sense that can be made precise, the helix-faced chiral polyhedra can be thought of as continuous ``chiral deformations" of helix-faced regular polyhedra (see \cite{pelwei}). This surprising phenomenon is illustrated for the helix-faced polyhedra $P_2(c,d)$ in Figure~\ref{figcont}; shown is the effect on the location of the ``vertical" helical faces, as a result of continuously changing the parameters $c,d$.

Finally, helix-faced chiral polyhedra are combinatorially regular, as they are isomorphic to regular polyhedra. However, by contrast, finite-faced chiral polyhedra are {\em combinatorially chiral\/}, meaning that the combinatorial automorphism group has two flag-orbits such that adjacent flags are in distinct orbits (see \cite{pelwei}).

In summary, we have the following

\begin{theorem}
\label{classifichiral}
Up to congruence, the chiral polyhedra in $\mathbb{E}^3$ fall into six infinite, $2$-parameter families of apeirohedra, each containing alongside chiral apeirohedra also two regular polyhedra. Three families consist of finite-faced apeirohedra, and three of helix-faced polyhedra. The finite-faced polyhedra are also combinatorially chiral, but the helix-faced polyhedra are combinatorially regular.
\end{theorem}

\section{Two-orbit polyhedra}
\label{2orbpoly}

The chiral polyhedra in $\mathbb{E}^3$ are by definition the $2$-orbit polyhedra in $\mathbb{E}^3$ in the class $2_I$ with $I=\emptyset$. It is desirable to extend the classification of chiral polyhedra to $2$-orbit polyhedra in arbitrary classes $2_I$, with $I\subsetneq\{0,1,2\}$. We saw that chirality cannot occur among finite polyhedra; however, as the example of the cuboctahedron (in class $2_{\{0,1\}}$) shows, finite $2$-orbit polyhedra already occur among the familiar convex polyhedra. Thus a good first step would be the complete enumeration of the finite $2$-orbit polyhedra in $\mathbb{E}^3$.

Significant progress towards this goal has already been made for regular polyhedra of index $2$. A polyhedron $\K$ is said to be a {\em regular polyhedron of index $2$\/} if its combinatorial automorphism group $\Gamma(\K)$ acts flag-transitively on $\K$ and contains the geometric symmetry group $G(\K)$ as a subgroup of index $2$. In other words, $\K$ is combinatorially regular but  ``fails geometric regularity by a factor of $2$". For any such polyhedron, the symmetry group has two orbits on the flags, and at most two orbits on the vertices, edges, and faces. Note that the helix-faced chiral polyhedra discussed in the previous section are examples of infinite regular polyhedra of index~$2$.

The finite regular polyhedra of index~$2$ were recently enumerated in Cutler-Schulte~\cite{cutsch} and Cutler~\cite{cut} (see also Wills~\cite{wills}). The following theorem summarizes the results.

\begin{theorem}
\label{classifindex2}
Up to similarity, there are exactly $22$ infinite families of regular polyhedra of index $2$ with vertices on two orbits under the symmetry group, where two polyhedra belong to the same family if they differ only in the relative size of the spheres containing their vertex orbits. In addition, up to similarity, there are exactly~$10$ (individual) polyhedra with vertices on one orbit under the symmetry group.
\end{theorem}

In describing the polyhedra, we slightly abuse terminology and say that a polyhedron $\K$ is of {\em type $\{p,q\}_r$} if the underlying regular map has (Schl\"afli) type $\{p,q\}$ and Petrie polygons of length $r$. Note here that we are not requiring the map to be the universal regular map of type $\{p,q\}$ with Petrie polygons of length $r$ (denoted $\{p,q\}_r$ in \cite{cm}). However, in some case the map actually is universal (see \cite{csw}).

\begin{figure}[htbp]
\begin{center}
\quad\includegraphics[width=1.75in]{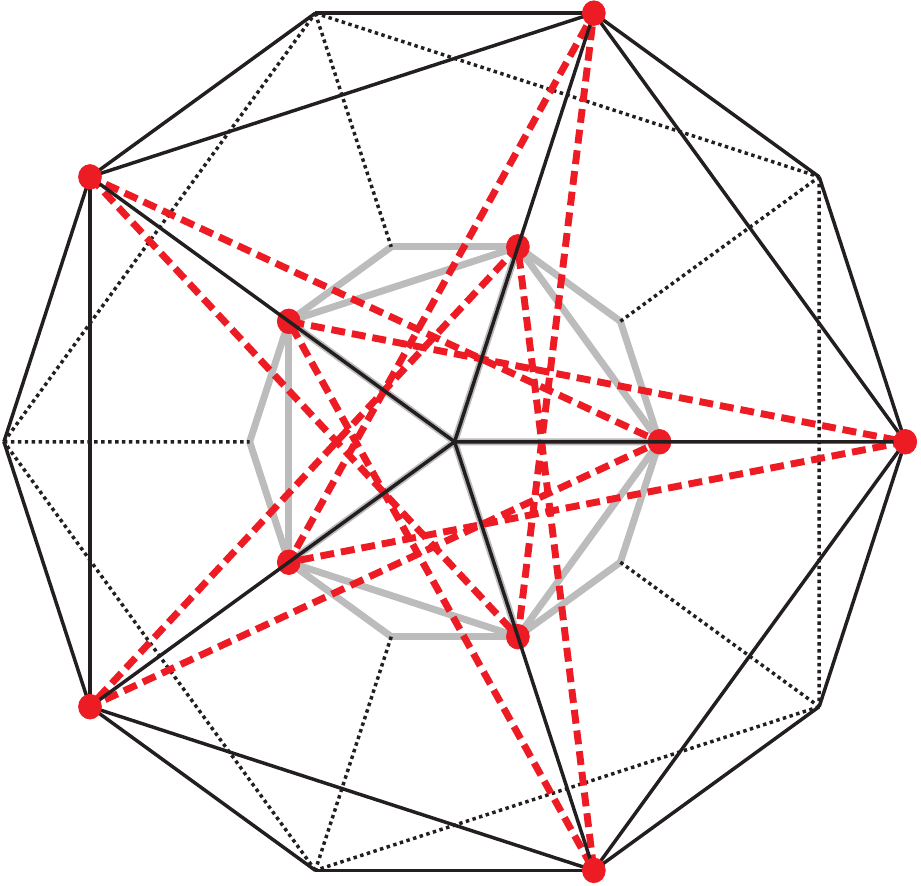}\qquad\qquad
\includegraphics[width=1.75in]{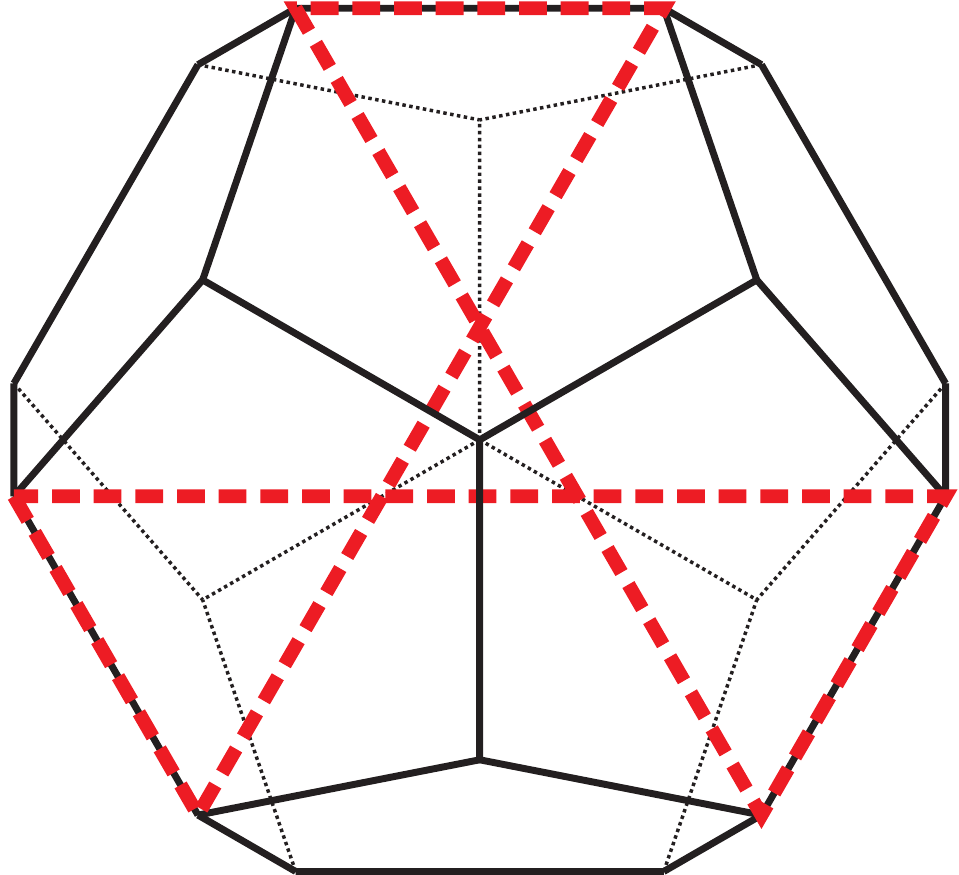}\qquad\quad
\end{center}
\caption{Two regular polyhedra of index $2$. The polyhedron on the left is a representative of one of the four infinite families of type $\{10,5\}_6$ with two vertex-orbits; its vertices lie on a pair of concentric icosahedra. The planar-faced polyhedron on the right has type $\{6,6\}_6$ and one vertex-orbit; its vertices are those of a dodecahedron. Only one face is shown in each case; the other faces are obtained by applying all icosahedral or dodecahedral symmetries.  }
\label{fig66}
\end{figure}

Table~\ref{tabfour} records the $22$ infinite families of polyhedra by combinatorial isomorphism type. For example, the last entry in row 5 indicates that there are $2$ infinite families with polyhedra isomorphic to Gordan's (universal) map $\{4,5\}_{6}$. The third column gives the name of the map in the notation of Conder~\cite{con} (when applicable), with $R$ or $N$, respectively, indicating an orientable or non-orientable regular map; the number before the period is the genus, and an asterisk indicates the dual. The polyhedra in these $22$ families have their vertices located at those of a pair of similar, aligned or opposed, Platonic solids with the same symmetry group. There are respectively $4$, $2$ and $16$ families with full tetrahedral, octahedral, and icosahedral symmetry. The symmetry group is face-transitive in each case, and each polyhedron is orientable. Among all polyhedra (in all families), there are just two polyhedra with planar faces. Figure~\ref{fig66} shows one face of a regular polyhedron of index $2$ and type $\{10,5\}_6$ belonging to one of the four families in the last row of Table~\ref{tabfour}.

\begin{table}[htbp]
\begin{center}
\begin{tabular}{|c|c|c|c|c|cl}
\hline
Type &Face Vector  &Map  & \# Families\\
$\{p,q\}_r$ &$(f_{0},f_{1},f_{2})$& &\\[.05in]
\hline
\hline
$\{4,3\}_6$   & $(8, 12, 6)$& sphere& 2\\
\hline
$\{6,3\}_4$   & $(8, 12, 4)$&torus &2\\
\hline
\hline
$\{6,4\}_6$   & $(12, 24, 8)$&$R3.4^*$ &2\\
\hline
\hline
$\{10,3\}_{10}$  & $(40,60,12)$&$R5.2^*$& 4\\
\hline
$\{4,5\}_{6}$  & $(24,60,30)$  &$R4.2$ &2\\
\hline
$\{6,5\}_{4}$   & $(24,60,20)$&$R9.16^*$ &2\\
\hline
$\{6,5\}_{10}$  & $(24,60,20)$&$R9.15^*$ & 4\\
\hline
$\{10,5\}_{6}$ & $(24,60,12)$ &$R13.8^*$ & 4\\
\hline
\end{tabular}
\end{center}
\caption{The $22$ infinite families of regular polyhedra of index $2$ with two vertex orbits, listed by combinatorial isomorphism type. The polyhedra in the first two rows have full tetrahedral symmetry, and those in the third row full octahedral symmetry; all others have full icosahedral symmetry.}
\label{tabfour}
\end{table}

The $10$ (individual) regular polyhedra of index $2$ with vertices on one orbit are listed in Table~\ref{tabfive}. Each has full icosahedral symmetry. There are orientable and non-orientable examples. Figure~\ref{fig66} depicts one face of the planar-faced regular polyhedron of index $2$ and type $\{6,6\}_6$ listed in the first row of Table~\ref{tabfive}.

\begin{table}[htb]
\begin{center}
\begin{tabular}{|c|c|c|c|}  \hline
Type & Face Vector  &Map& Notes \\
$\{p,q\}_r$ & $(f_{0},f_{1},f_{2})$&&\\[.05in]
\hline
\hline
$\{6, 6\}_6$ & $(20, 60, 20)$&$R11.5$& planar faces,\\
&&&self-dual map\\
\hline
$\{6, 6\}_6$& $(20, 60, 20)$ &$N22.3$ & face transitive\\
\hline
$\{4, 6\}_5$& $(20, 60, 30)$& $N12.1$ &\\
\hline
$\{5, 6\}_4$& $(20, 60, 24)$& $R9.16$& planar faces\\
\hline
$\{6, 4\}_5$& $(30, 60, 20)$& ${N12.1}^*$ &\\
\hline
$\{5, 4\}_6$& $(30, 60, 24)$& $R4.2^*$& planar faces\\
\hline
$\{4, 6\}_{10}$&$(20, 60, 30)$& $R6.2$ &\\
\hline
$\{10, 6\}_4$& $(20, 60, 12)$& $N30.11^*$& \\
\hline
$\{6, 4\}_{10}$&$(30, 60, 20)$& $R6.2^{*}$&\\
\hline
$\{10, 4\}_6$& $(30, 60, 12)$& $N20.1^*$&\\
\hline
\end{tabular}
\end{center}
\caption{The $10$ (individual) regular polyhedra of index $2$ with one vertex orbits. Each has full icosahedral symmetry.}
\label{tabfive}
\end{table}

\section{Conclusions}

The recent history of symmetric structures in Euclidean 3-space $\mathbb{E}^3$ suggests a rich variety of objects yet to be discovered. All geometrically regular polygonal complexes (including polyhedra), and all regular 4-polytopes, in Euclidean 3-space have now been classified; by contrast, little is known about polygonal complexes and 4-polytopes with slightly less symmetry. Two natural open questions concern the enumeration of all 2-orbit polyhedra and all edge-transitive polyhedra in $\mathbb{E}^3$. It appears more challenging to widen the scope of these problems to general polygonal complexes. A good starting point in this direction is a detailed classification of the finite 2-orbit, or edge-transitive, polygonal complexes that are not polyhedra; or a proof that such complexes cannot exist.

Significant progress has been made in the theory of realizations (in any dimension) for regular polytopes of any rank, mostly by McMullen; the state of the art will be summarized in his forthcoming monograph on ``Geometric Regular Polytopes"~\cite{grp}, but many results can also be found in \cite[Chapter 5]{arp}. However, little is known about realizations of other kinds of polytopes or polygonal complexes. The complete enumeration of particularly interesting families of such objects will greatly contribute to our basic understanding of geometric realizations of incidence structures.

\begin{acknowledgement}
We greatly appreciate the hospitality of the Fields Institute over extended periods of time during the Thematic Program on Discrete Geometry and Application in Fall 2011, and are very grateful for the support we have received. Daniel Pellicer was a postdoctoral fellow at Fields Institute in Fall 2011, and was also partially supported by PAPIIT--Mexico under grant IN106811-3 and CONACYT project 166951. Egon Schulte was also supported by NSF-grant DMS--0856675.
\end{acknowledgement}


\begin{thebibliography}{99.}
\bibitem{ar} J.L.Arocha, J.Bracho and L.Montejano, {\em Regular projective polyhedra with planar faces, Part I\/}, Aequat.\ Math.\  59 (2000), 55--73.

\bibitem{bra} J.Bracho, {\em Regular projective polyhedra with planar faces, Part II\/}, Aequat.\ Math.\ 59 (2000), 160--176.

\bibitem{bu2} F.Buekenhout (ed.), {\em Handbook of Incidence Geometry\/}, Elsevier Science B.V. (1995).

\bibitem{con} M.~Conder, {\em Regular maps and hypermaps of Euler characteristic $-$1 to $-$200}, J. Comb. Theory Ser.~B 99 (2009), 455--459.  (Associated lists available online: http://www.math.auckland.ac.nz/$\sim$conder.)

\bibitem{crsp} H.S.M.Coxeter, {\em Regular skew polyhedra in 3 and 4 dimensions and their topological analogues\/},  Proc.\ London Math.\ Soc. (2) 43 (1937), 33--62.  (Reprinted with amendments in {\em Twelve Geometric Essays\/}, Southern Illinois University Press (Carbondale, 1968), 76--105.)

\bibitem{coxeter} H.S.M.Coxeter, {\em Regular Polytopes\/} (3rd edition),  Dover
(New York, 1973).

\bibitem{cm} H.S.M.Coxeter and W.O.J.Moser, {\em Generators and Relations for Discrete Groups}, 4th editon, Springer, Berlin, 1980.

\bibitem{cut} A.M.Cutler, {\em Regular polyhedra of index two, II\/}, Contributions to Algebra and Geometry 52 (2011), 357--387.

\bibitem{cutsch} A.M.Cutler and E.Schulte, {\em Regular polyhedra of index two, I\/}, Contributions to Algebra and Geometry 52 (2011), 133--161.

\bibitem{csw} A.M.Cutler, E.Schulte and J.M.Wills, {\em Icosahedral skeletal polyhedra realizing Petrie relatives of Gordan's regular map\/}, Contributions to Algebra and Geometry (to appear).

\bibitem{kom1} L.Danzer and E.Schulte, {\em Regul\"are Inzidenzkomplexe, I\/},
Geom.\ Dedicata 13 (1982), 295--308.

\bibitem{del} O.Delgado-Friedrichs, M.D.Foster, M.O'Keefe, D.M.Proserpio, M.M.J.Treacy and O.M.Yaghi, {\em What do we know about three-periodic nets?\/}, J. Solid State Chemistry 178 (2005), 2533--2554.

\bibitem{d1} A.W.M.Dress, {\em A combinatorial theory of Gr\"unbaum's new regular polyhedra, I:  Gr\"unbaum's new regular polyhedra and their automorphism group\/},  Aequationes Math.\  23 (1981), 252--265.

\bibitem{d2} A.W.M.Dress, {\em A combinatorial theory of Gr\"unbaum's new regular polyhedra, II:  complete enumeration\/},  Aequationes Math.\ 29 (1985), 222--243.

\bibitem{grove} L.C.Grove and C.T.Benson, {\em Finite Reflection Groups\/}
(2nd edition),  Graduate Texts in Mathematics, Springer-Verlag (New
York-Heidelberg-Berlin-Tokyo, 1985).

\bibitem{gr1} B.Gr\"unbaum, {\em Regular polyhedra --- old and new\/}, Aequat.\ Math.\ {16} (1977), 1--20.

\bibitem{grhol} B.Gr\"unbaum, {\em Polyhedra with hollow faces\/},  in {\em
Polytopes:  Abstract, Convex and Computational\/} (eds.\ T.~Bisztriczky,
P.~McMullen, R.~Schneider and A.~Ivi\'c Weiss), NATO ASI Series C 440,
Kluwer (Dordrecht etc., 1994), 43--70.

\bibitem{gr} B.Gr\"unbaum, {\em Acoptic polyhedra\/}, In {\em Advances in Discrete and
Computational Geometry\/}, B.Chazelle et al. (ed.), Contemp.\ Math.\ 223,
American Mathematical Society (Providence, RI, 1999), 163--199.

\bibitem{hub} I.Hubard, {\em Two-orbit polyhedra from groups\/}, European Journal of Combinatorics.
31 (2010), 943--960.

\bibitem{hubsch} I.Hubard and E.Schulte, {\em Two-orbit polytopes\/}, in preparation.

\bibitem{pm} P.McMullen, {\em Regular polytopes of full rank\/}, Discrete \& Computational Geometry 32 (2004), 1--35.

\bibitem{grp} P.McMullen, {\em Geometric Regular Polytopes\/}, in preparation.

\bibitem{ordinary} P.McMullen and E.Schulte, {\em Regular polytopes in ordinary space\/}, Discrete Comput. Geom. 17 (1997), 449--478.

\bibitem{arp}
P.McMullen and E. Schulte, \emph{Abstract regular polytopes}, Encyclopedia of Mathematics and its Applications, Vol. 92, Cambridge University Press, Cambridge, UK, 2002.

\bibitem{ms3} P.McMullen and E.Schulte, {\em Regular and chiral polytopes in low dimensions\/}, In {\em The Coxeter Legacy -- Reflections and Projections\/} (eds. C.Davis and E.W.Ellers), Fields Institute Communications, Volume 48, American Mathematical Society (Providence, RI, 2006), 87--106.

\bibitem{monw5} B.R.Monson and A.I.Weiss, {\em Realizations of regular toroidal maps\/},
Canad.\ J.\ Math. (6) 51 (1999), 1240--1257.

\bibitem{okee}
M.O'Keeffe, {\em Three-periodic nets and tilings:\  regular and related infinite polyhedra\/}, Acta Crystallographica  A 64 (2008), 425--429.

\bibitem{okhy} M.O'Keeffe and B.G.Hyde, {\em Crystal Structures; I. Patterns and Symmetry\/}, Mineralogical Society of America, Monograph Series, Washington, DC, 1996.

\bibitem{pelsch1} D.Pellicer and E.Schulte, {\em Regular polygonal complexes in space, I}, Trans. Amer. Math. Soc. 362 (2010), 6679--6714.

\bibitem{pelsch2} D.Pellicer and E.Schulte, {\em Regular polygonal complexes in space, II}, Trans. Amer. Math. Soc. (to appear).

\bibitem{pelwei} D.Pellicer and A.I.Weiss, {\em Combinatorial structure of Schulte's chiral polyhedra\/}, Discrete \& Comput. Geom. 44 (2010), 167--194.

\bibitem{ratcliffe}
J.G.Ratcliffe, {\em Foundations of Hyperbolic Manifolds\/},
Graduate Texts in Mathematics, Springer-Verlag (New York-Berlin-Heidelberg, 1994).

\bibitem{sympopo}  E.Schulte, {\em Symmetry of polytopes and polyhedra\/},  In {\em Handbook of Discrete and Computational Geometry, Second Edition\/} (eds. J.E.Goodman and J.O'Rourke), Chapman \& Hall/CRC (Boca Raton, 2004), 431--454.

\bibitem{kom2}
E.Schulte, {\em Regul\"are Inzidenzkomplexe, II\/},
Geom.\ Dedicata 14 (1983), 33--56.

\bibitem{chiral1}
E.Schulte, {\em Chiral polyhedra in ordinary space, I\/}, Discrete Comput. Geom. 32
(2004), 55--99.

\bibitem{chiral2}
E.Schulte, {\em Chiral polyhedra in ordinary space, II\/}, Discrete Comput. Geom. 34
(2005), 181--229.

\bibitem{wells} A.F.Wells, {\em Three-dimensional Nets and Polyhedra\/},
Wiley-Interscience (New York, etc., 1977).

\bibitem{wills} J.M.Wills, {\em Combinatorially regular polyhedra of index~2\/}, Aequat. Math. 34 (1987), 206--220.

\end{thebibliography}
\end{document}